\documentclass{article}
\usepackage{amssymb}
\usepackage{amsmath}
\usepackage{amsthm}
\usepackage{mathrsfs}
\usepackage{mathtools}
\usepackage{tikz-cd}
\usepackage{quiver}

\usepackage[hidelinks]{hyperref}

\theoremstyle{definition}
\newtheorem{theorem}{Theorem}

\newtheorem{definition}{Definition}
\newtheorem*{definition*}{Definition}
\theoremstyle{remark}
\newtheorem*{example}{Example}

\newcommand{\Cat}[1]{{\mathcal{#1}}}
\DeclareMathOperator{\Hom}{Hom}

\title{Locally Cartesian Closed Categories}
\author{Huang Xu\footnote{Many thanks to Alias Qli who
helped proofread and revise this.}}
\begin{document}
\maketitle

\begin{abstract}
This note explains how dependent sums and products
are interpreted by adjoints of the base change functor
in a locally cartesian closed category. An effort is made
to unpack all the definitions so as to make the concepts
more transparent to new learners.
\end{abstract}

Notational conventions:
\begin{itemize}
    \item Categories in general use the calligraphic font: \(\Cat C\);
    \item Special categories use sans-serif: \(\mathsf{Cat}, \mathsf{Set}\);
    \item \(X \in \Cat C\) means ``\(X\) is an object in \(\Cat C\)\,'';
    \item Compositions are in the ``function order'', i.e. if \(X \xrightarrow{f} Y \xrightarrow{g} Z\),
    then the composite arrow is \(g\circ f\);
    \item \(\Hom(X, Y)\) denotes the morphisms from \(X\) to \(Y\);
    If necessary, subscripts indicate the category in discussion: \(\Hom_{\Cat C}(X, Y)\).
\end{itemize}

\section{Slices}
Given a category \(\Cat C\) and an object \(X\),
let's consider all the arrows into \(X\). This forms a
collection of arrows \[\bigcup_{Y\in \Cat C}\Hom(Y, X).\]
We shall take this collection of arrows as the \emph{objects}
of a new category, named \(\Cat C/X\).

What should the morphisms be? Consider any commutative
diagram of the form

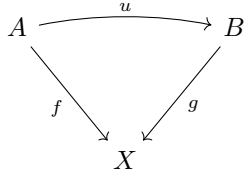
\[\begin{tikzcd}
    A && B \\
    \\
    & X
    \arrow["f"', from=1-1, to=3-2]
    \arrow["g", from=1-3, to=3-2]
    \arrow["u", curve={height=-6pt}, from=1-1, to=1-3]
\end{tikzcd}\]

By saying that the diagram ``commutes'', I mean
\(g \circ u = f\). It would be natural to take \(u\)
as a morphism from \(f\) to \(g\). This defines the
\textbf{slice} category over \(X\).

\begin{example}
Here are some simple examples of slice categories.
\begin{itemize}
\item If \(\Cat C\) has a terminal object, then \(\Cat C/1 \cong \Cat C\).
\item Take \(2\) to be the set \(\{\mathsf{blue}, \mathsf{red}\}\).
\(\mathsf{Set}/2\) is the category of \emph{two-colored sets}.
In other words, its objects are sets where each element is assigned
either the color \(\mathsf{blue}\) or \(\mathsf{red}\).
Morphisms are set-theoretic functions that maps blue elements
to blue ones, and vice versa.
\item \(\mathsf{Set}/\varnothing\) contains only one object and one morphism.
\item \textsc{Exercise}: Come up with one more example. Make it as interesting as you can.
\end{itemize}
\end{example}

Notice that given \emph{any} object in a category, we can
make a slice category out of it. So suppose we have two
objects and a morphism \(X \xrightarrow{f}Y\). What can we say
of the two slice categories?

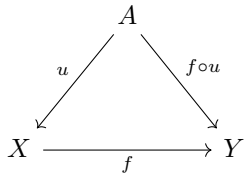
\[\begin{tikzcd}
    & A \\
    \\
    X && Y
    \arrow["u"', from=1-2, to=3-1]
    \arrow["f\circ u", from=1-2, to=3-3]
    \arrow["f"', from=3-1, to=3-3]
\end{tikzcd}\]

Here, \(u \in \Cat C/X\) and \(f\circ u \in \Cat C/Y\).
Therefore, there is a map from the objects of \(\Cat C/X\)
to the objects of \(\Cat C/Y\). The next question to ask,
is whether the map is \emph{functorial}.
Here's the relevant diagram. The verification is left as an exercise.

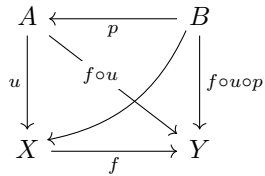
\[\begin{tikzcd}
    A && B \\
    \\
    X && Y
    \arrow["u"', from=1-1, to=3-1]
    \arrow["{f\circ u}"{description, pos=0.4}, from=1-1, to=3-3]
    \arrow["f"', from=3-1, to=3-3]
    \arrow["p", from=1-3, to=1-1]
    \arrow["{f\circ u\circ p}", from=1-3, to=3-3]
    \arrow[shift right=1, curve={height=-12pt}, from=1-3, to=3-1]
\end{tikzcd}\]

We give this functor a name: \(f_! : \Cat C/X \to \Cat C/Y\).

\section{Pullbacks}

The next thing we do requires more structure in the category \(\Cat C\).
Let's take three objects \(B \xrightarrow{f} A \xleftarrow{g} C\).
If there happens to exist \(X\) together with arrows
\(B \xleftarrow{p} X \xrightarrow{q} C\) such that the square commutes,
and additionally...

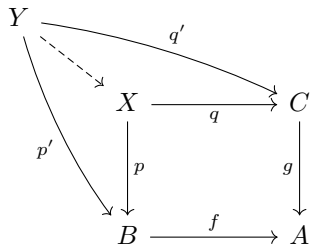
\[\begin{tikzcd}
    Y \\
    & X && C \\
    \\
    & B && A
    \arrow["p", from=2-2, to=4-2]
    \arrow["f", from=4-2, to=4-4]
    \arrow["q"', from=2-2, to=2-4]
    \arrow["g"', from=2-4, to=4-4]
    \arrow[dashed, from=1-1, to=2-2]
    \arrow["{p'}"', curve={height=6pt}, from=1-1, to=4-2]
    \arrow["{q'}", curve={height=-6pt}, from=1-1, to=2-4]
\end{tikzcd}\]

... For every given \(Y\) with morphisms \(p',q'\), there is a
unique arrow \(Y \to X\) such that the diagram commutes. In
this case, we call \(X\) a \textbf{pullback}.

What are pullbacks like? We need to find arrows
\(p, q\) that ``reconcile'' \(f\) and \(g\). In \(\mathsf{Set}\),
the pullback is given by the set
\[\{(b, c) \mid f(b) = g(c)\},\]
equipped with the obvious projections \(p, q\).

But there's another way to look at it. Each point \(a \in A\)
determines a set \(f^{-1}(a) = \{ b \mid f(b) = a\}\), and similarly
\(g^{-1}(a)\). This is called the \textbf{preimage}.
In this way, \(B\) can be rewritten as a union of preimages:
\[ B = \bigcup_{a \in A} f^{-1}(a). \]
\textsc{Exercise}: In this union, each set is disjoint from each other. Can
you see why? Since they are disjoint, we can use \(\coprod\) instead
of \(\bigcup\) to emphasize this (these two symbols have the same
meaning except \(\coprod\) implies disjointness).

Therefore, we may regard \(B\) as a space composed of ``fibers''
\(f^{-1}(a)\). For example, if \(B = \mathbb R^2\), and \(A = \mathbb R\),
take \[f(x, y) = x^2 + y^2.\] Then \(B\) is divided into concentric
circles \(f^{-1}(r^2)\) of radius \(r\) about the origin.
Note that \(f^{-1}(-1)\) is empty, meaning that the fiber that
lies over \(-1\) is \(\varnothing\).

What does this has to do with pullbacks? Well, we can rewrite
\(X\) in this way:
\[X = \coprod_{a\in A} f^{-1}(a) \times g^{-1}(a).\]
It is another fibered space, where each fiber is the \emph{product}
of the corresponding fibers in \(B\) and \(C\).
From this perspective, we may call the pullback as \textbf{fibered product},
denoted \(B \times_A C\). \textsc{Exercise}: Prove that \(B \times_1 C \cong B \times C\)
holds in any category with a terminal object.

The reader should be familiar with the fact that
\(A \times (-)\) is a functor. This is in accordance with
the Haskell typeclass instance \texttt{Functor ((,) a)}.
In fact, pullbacks, being called the fibered product, is also
a functor. To verify this, we need a diagram:

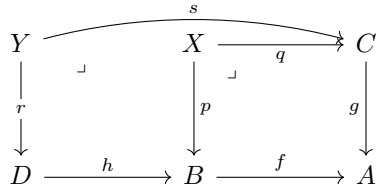
\[\begin{tikzcd}
    Y && X && C \\
    \\
    D && B && A
    \arrow["p", from=1-3, to=3-3]
    \arrow["f", from=3-3, to=3-5]
    \arrow["q"', from=1-3, to=1-5]
    \arrow["g"', from=1-5, to=3-5]
    \arrow["h", from=3-1, to=3-3]
    \arrow["s", curve={height=-12pt}, from=1-1, to=1-5]
    \arrow["r"{description}, from=1-1, to=3-1]
    \arrow["\lrcorner"{anchor=center, pos=0.125}, draw=none, from=1-3, to=3-5]
    \arrow["\lrcorner"{anchor=center, pos=0.125}, draw=none, from=1-1, to=3-5]
\end{tikzcd}\]

Here the little right-angle marks say that there are two
pullback squares. We need to prove that there is
an arrow \((\mathsf{fmap}\, h) : Y \to X\). This follows
directly from the universal property of pullbacks. Next,
we need the functor law.

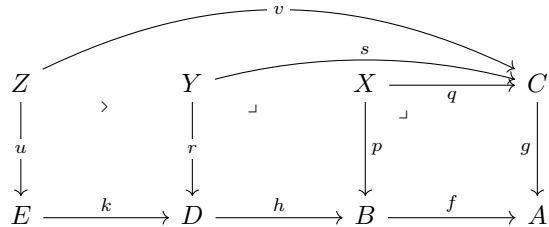
\[\begin{tikzcd}
    Z && Y && X && C \\
    \\
    E && D && B && A
    \arrow["p", from=1-5, to=3-5]
    \arrow["f", from=3-5, to=3-7]
    \arrow["q"', from=1-5, to=1-7]
    \arrow["g"', from=1-7, to=3-7]
    \arrow["h", from=3-3, to=3-5]
    \arrow["s", curve={height=-12pt}, from=1-3, to=1-7]
    \arrow["r"{description}, from=1-3, to=3-3]
    \arrow["\lrcorner"{anchor=center, pos=0.125}, draw=none, from=1-5, to=3-7]
    \arrow["\lrcorner"{anchor=center, pos=0.125}, draw=none, from=1-3, to=3-7]
    \arrow["k", from=3-1, to=3-3]
    \arrow["u"{description}, from=1-1, to=3-1]
    \arrow["v"{description}, shift left=1, curve={height=-30pt}, from=1-1, to=1-7]
    \arrow["\lrcorner"{anchor=center, pos=0.125, rotate=45}, draw=none, from=1-1, to=3-7]
\end{tikzcd}\]

The reader shall complete the argument using the given diagram.

Before we move on, let's pause for a moment and ponder what
we just proved. Note that to use \((\mathsf{fmap}\,h)\),
the large square \(Y,C,D,A\) cannot be arbitrary:
The lower edge has to be \(D\xrightarrow{f\circ h}A\).
So what is the ``functor'' that we've just found?
What is its source and target categories? It turns out that
\((-) \times_A C\) is actually a functor \(\Cat C/A \to \Cat C/C\)!
The choice of these categories are important.
Note that although in the notation \(B \times_A C\), the two
arrows \(f, g\) doesn't appear, they are the essential
ingredients. \textsc{Exercise}: Give an example
of two pullbacks \(B \times_A C\) with different \(g\), such
that the results are not isomorphic.

Saying that the functor is in \(\Cat C/A \to \Cat C/C\)
instead of \(\Cat C \to \Cat C\) adds the important information
of the respective arrows into \(A\). And this ensures that
a morphism in \(\Cat C/A\) always commutes with these arrows.

To emphasize the importance of the morphisms, we write \(g^* :\Cat C/A \to \Cat C/C\)
for the functor. Note that the functor goes in the opposite
direction of \(g : C \to A\). But this does \emph{not} make
\(g^*\) a contravariant functor. As you have proved in the previous
section, \(g^*\) turns \(h : M \to N\)
into \((\mathsf{fmap}\,g^*)h : M \times_A C \to N \times_A C\),
which means it is covariant.

\section{Adjoint Yoga}

Anyway, we now have a functor \(g_! : \Cat C/C \to \Cat C/A\)
from the first section, and \(g^* : \Cat C/A \to \Cat C/C\)
from the second section. In category theory, whenever you encounter
this, make a bet that they are adjoint.

What is adjunction? There are two equivalent definitions
that I find the most natural. The first one describes an
adjoint pair as an \emph{almost} inverse pair of functors.

\begin{definition}
Two functors \(F : \Cat C \to \Cat D\) and
\(G : \Cat D \to \Cat C\) are called \textbf{adjoint} if the following holds.

For each object \(X \in \Cat D\),
there is a morphism \(\epsilon_X : FGX \to X\), and similarly
for each \(Y \in \Cat C\) a \(\eta_Y : Y \to GFY\), satisfying
the following conditions:
\begin{itemize}
\item The assignment of morphisms \(\epsilon_X\) is natural.
In other words, for a morphism \(f : X_1 \to X_2\), we have
\(\mathsf{fmap}_{FG}f : FGX_1 \to FGX_2\), this forms a square
\[\begin{tikzcd}
    {FGX_1} && {FGX_2} \\
    \\
    {X_1} && {X_2}
    \arrow["{\epsilon_{X_1}}"{description}, from=1-1, to=3-1]
    \arrow["{\epsilon_{X_2}}"{description}, from=1-3, to=3-3]
    \arrow["{\mathsf{fmap}\, f}"{description}, from=1-1, to=1-3]
    \arrow["f"{description}, from=3-1, to=3-3]
\end{tikzcd}\]
The naturality condition states that all these squares commute.
Similar conditions hold for \(\eta\).
\item \(\epsilon,\eta\) settles the situation for composing \emph{two}
functors. In the case of three functors, we have two maps
\[FGFX \xleftrightharpoons[\epsilon_{FX}]{\mathsf{fmap}\, \eta_{X}} FX.\]
These should compose to get the identity on \(FX\). Similar conditions
hold fo \(GY\).
\end{itemize}

In this case, \(F\) is called the left adjoint, and \(G\) the right adjoint,
denoted as \(F \dashv G\).
\end{definition}

I won't linger too much on the concept of adjunction. But here're
two quick examples.
\begin{itemize}
\item \(U : \mathsf{Mon} \to \mathsf{Set}\) is
a functor that maps a monoid to its underlying set. And
\(F : \mathsf{Set} \to \mathsf{Mon}\) maps a set \(X\)
to the collection of lists \(\mathsf{List}(X)\), with list
concatenation as monoid multiplication, and the
empty list \([]\) as the neutral element. \(F\) is left adjoint
to \(U\).
\item Let \(\Delta : \mathsf{Set} \to \mathsf{Set}\times\mathsf{Set}\)
be the diagonal functor, sending \(X\) to \((X,X)\).
The product functor \((-)\times(-) : \mathsf{Set}\times\mathsf{Set} \to \mathsf{Set}\)
is the right adjoint of \(\Delta\).
\end{itemize}

The second definition is more catchy:

\begin{definition}
Two functors \(F : \Cat C \to \Cat D\) and
\(G : \Cat D \to \Cat C\) are adjoint iff
\[\Hom(FX, Y) \cong \Hom(X,GY)\]
such that the isomorphism is natural in \(X\) and \(Y\).
\end{definition}

The reader shall verify that these two definitions are equivalent,
and that the two examples given are indeed adjoints (using both definitions).

Now let's turn back to our two functors \(f_!, f^*\).
We draw a diagram to compose them and see what happens.
First look at \(f_!f^*x\).

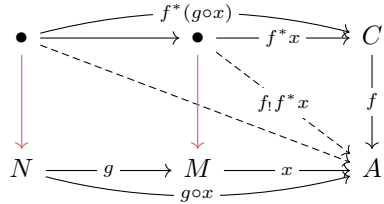
\[\begin{tikzcd}
    \bullet && \bullet && C \\
    \\
    N && M && A
    \arrow["f"{description}, from=1-5, to=3-5]
    \arrow["x"{description}, from=3-3, to=3-5]
    \arrow["{g\circ x}"{description}, curve={height=12pt}, from=3-1, to=3-5]
    \arrow["g"{description}, from=3-1, to=3-3]
    \arrow["{f^*x}"{description}, from=1-3, to=1-5]
    \arrow[from=1-1, to=1-3]
    \arrow["{f^*(g\circ x)}"{description}, curve={height=-12pt}, from=1-1, to=1-5]
    \arrow[dashed, from=1-1, to=3-5]
    \arrow["{f_!f^*x}"{description}, dashed, from=1-3, to=3-5]
    \arrow[color={rgb,255:red,214;green,92;blue,92}, from=1-1, to=3-1]
    \arrow[color={rgb,255:red,214;green,92;blue,92}, from=1-3, to=3-3]
\end{tikzcd}\]

The lower half is in \(\Cat C/A\), and the upper half in
\(\Cat C/C\). The two dashed arrows are \(x\) and \(g\circ x\)
under the functor \(f_!f^*\). They lie in \(\Cat C/A\).
Now notice the red arrows generated from the pullback.
Composing them with each \(x \in \Cat C/A\)
gives a transformation from \(x\) to \(f_!f^*x\).
This gives \(\eta_x : x \to f_!f^*x\).

What about the naturality condition? \textsc{Exercise}: Argue that
the square \(\bullet,\bullet,N,M\) commutes, and explain
why this proves the naturality condition for \(\eta\).

Next, the reverse composition \(f^*f_!\). It is slightly
trickier:

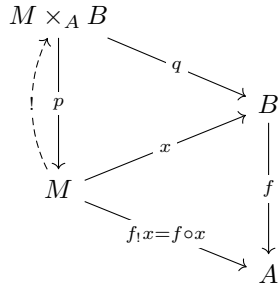
\[\begin{tikzcd}
    {M\times_AB} \\
    && B \\
    M \\
    && A
    \arrow["f"{description}, from=2-3, to=4-3]
    \arrow["x"{description}, from=3-1, to=2-3]
    \arrow["{f_!x=f\circ x}"{description}, from=3-1, to=4-3]
    \arrow["p"{description}, from=1-1, to=3-1]
    \arrow["q"{description}, from=1-1, to=2-3]
    \arrow["{!}"{description}, curve={height=-12pt}, dashed, from=3-1, to=1-1]
\end{tikzcd}\]

Here we have \(x \in \Cat C/B\). Therefore, there is
a well-hidden commutative square:

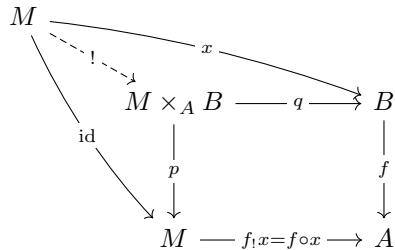
\[\begin{tikzcd}
    M \\
    & {M\times_AB} && B \\
    \\
    & M && A
    \arrow["f"{description}, from=2-4, to=4-4]
    \arrow["{f_!x=f\circ x}"{description}, from=4-2, to=4-4]
    \arrow["p"{description}, from=2-2, to=4-2]
    \arrow["q"{description}, from=2-2, to=2-4]
    \arrow["{\mathrm{id}}"{description}, curve={height=6pt}, from=1-1, to=4-2]
    \arrow["x"{description}, curve={height=-6pt}, from=1-1, to=2-4]
    \arrow["{!}"{description}, dashed, from=1-1, to=2-2]
\end{tikzcd}\]

... Which creates the unique morphism \(!\), such that
\(p \circ {!} = \mathrm{id}\) and \(q \circ {!} = x\).
Now recall that \(q = f^*f_!x\). Therefore, composing
with \(!\) gives a natural transformation
\(\epsilon_x : f^*f_!x \to x\).

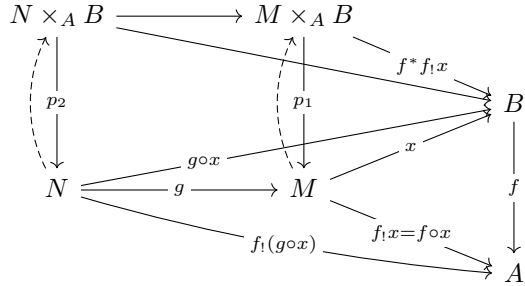
\[\begin{tikzcd}
    {N\times_A B} && {M\times_AB} \\
    &&&& B \\
    N && M \\
    &&&& A
    \arrow["f"{description}, from=2-5, to=4-5]
    \arrow["x"{description}, from=3-3, to=2-5]
    \arrow["{f_!x=f\circ x}"{description}, from=3-3, to=4-5]
    \arrow["{p_1}"{description}, from=1-3, to=3-3]
    \arrow["{f^*f_!x}"{description}, from=1-3, to=2-5]
    \arrow[curve={height=-12pt}, dashed, from=3-3, to=1-3]
    \arrow["g"{description}, from=3-1, to=3-3]
    \arrow["{g\circ x}"{description, pos=0.3}, from=3-1, to=2-5]
    \arrow["{f_!(g\circ x)}"{description}, curve={height=6pt}, from=3-1, to=4-5]
    \arrow[from=1-1, to=2-5]
    \arrow["{p_2}"{description}, from=1-1, to=3-1]
    \arrow[from=1-1, to=1-3]
    \arrow[curve={height=-12pt}, dashed, from=3-1, to=1-1]
\end{tikzcd}\]

The naturality condition amounts
to proving that the two dashed arrows form a commutative
square. This follows immediately from the universal
property of pullbacks.

If you find this dizzying, why not try the other definition?

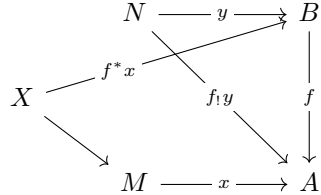
\[\begin{tikzcd}
    & N && B \\
    X \\
    & M && A
    \arrow["x"{description}, from=3-2, to=3-4]
    \arrow["y"{description}, from=1-2, to=1-4]
    \arrow["f"{description}, from=1-4, to=3-4]
    \arrow[from=2-1, to=3-2]
    \arrow["{f^*x}"{description, pos=0.3}, from=2-1, to=1-4]
    \arrow["{f_!y}"{description}, from=1-2, to=3-4]
\end{tikzcd}\]

You need to find a natural isomorphism between
\(\{g \mid y = f^* x \circ g\}\) and \(\{g \mid g\circ x = f\circ y\}\).
One direction is given by composition, and the other is given
by the universal property of pullbacks.

\section{Dependent Sum}

It's time to reveal the meaning of these constructions.
Recall how we can regard a morphism \(p : E \to B\)
as a fibered space \[E = \coprod_{x : B} p^{-1}(x).\]
So in the slice category \(\Cat C/B\), everything is fibered
along \(B\). If we take the map \({!}:B\to1\), then it
induces the functor \(\Cat C/B \to \Cat C/1\).
which takes a fibered space \(p : E \to B\) to \(E \to 1\).

Although this looks trivial, looking from the perspective
of fibered spaces, we get something different: \(p : E \to B\)
describes \(E\) with fibers over \(B\). And the functor
turns it into \(! : E \to 1\), where all the fibers are
merged into one big component. This corresponds to the
\textbf{dependent sum}:

\[ \sum_{x:B}p^{-1}(x). \]

We can generalize this by replacing the
terminal object with an arbitrary object \(A\),
and the morphism \({!} : B \to 1\) with
an arbitrary morphism
\(f : B \to A\), whose induced functor \(f_!\)
takes a ``fiberwise dependent sum'', i.e. for each
\(a \in A\), the fiber over \(a\) is
\[\sum_{x:B_a} p^{-1}(x),\]
where \(B_a\) is the fiber of \(B\) over \(a\).

What, then, is the functor \(f^*\)? Similarly we first take
\(A = 1\), and let \(f\) be
the unique morphism \({!} : B \to 1\). The pullback functor takes
\(p' : E \to 1\)
to \(\pi_1 : B \times E \to B\) projecting
to the first component.\footnote{Note that now
\(p' \in \Cat C/A\) (and we are studying the special
case \(A = 1\)), where in the last paragraph
\(p \in \Cat C/B\). This is because the functor
\(f^*\) goes in the opposite direction of \(f_!\),
and we need \(p'\) to be in the \emph{source} category
of the functor we are discussing.}
In the fibered space language,
it creates a trivial fibered space where each fiber looks
identical to \(E\).

Now generalizing to arbitrary \(f : B \to A\),
the pullback functor takes \(p' : E \to A\)
to a morphism \(E \times_A B \to B\). In the category
\(\mathsf{Set}\), the fibers of the new space looks
like \[{p'}^{-1}(f(b))\] for each \(b \in B\). In effect,
it changes the \emph{base space} from \(A\) to \(B\).
And thus it is named the \textbf{base change} functor.

\section{Towards Dependent Product}

The next goal is to characterize dependent products.
Following our previous experiences, it should
be a functor \(f_* : \Cat C/B \to \Cat C/A\) for \(f : B \to A\).
Similar to the dependent sum functor, it should
take a ``fiberwise dependent product'':
\[\prod_{x:B_a}p^{-1}(x),\]
where \(p : E \to B\) is regarded as a fibered space over \(B\).
As usual, we should consider the easy case where \(A = 1\),
and we only need to construct
\[\prod_{x:B}p^{-1}(x).\]

How should it be defined? \(\prod_{x:B}M\), where
\(M\) does not depend on \(x\),
is exactly the function space \(M^B\).
This suggests that we can define the dependent
product set \(\prod_{x:B}p^{-1}(x)\) as a subset of
the functions \(B \to \coprod_{x:B}p^{-1}(x)\). Of course, to
be type-correct, it needs to map \(b\in B\) to an
element of \(p^{-1}(b)\).
This can be expressed as it being a right inverse of
\(p\). So to sum up,
our quest is now to find right inverses \({?} \circ p = \mathrm{id}\)
of \(p\).

\subsection*{Interlude: Exponentials}

Actually, we not only need to find the right inverses.
In \(\mathsf{Set}\), we need a \emph{set} of right inverses,
which means instead of a collection of morphisms we need a
\emph{single object} that stands for the set of right inverses.
Before we tackle that, we shall look at how we can create a
single object that stands for the set of functions
--- the \emph{exponential object}.

How should a set of functions behave? Given sets \(X, Y\),
if we have a set of functions \(E = Y^X\), then we should
be able to \emph{evaluate} the functions at a given point
\(x\in X\). This is called the \emph{evaluation functional}%
\footnote{The ``-al'' part of the word ``functional''
is just something that stuck with mathematicians. It doesn't
really mean anything special.}
\[\mathrm{ev}(-,-) : E \times X \to Y.\]
So we already have the first parts of the definition:

\begin{definition}
Given objects \(X, Y\), an \textbf{exponential object} is defined
as an object \(E\) equipped with a morphism \(\mathrm{ev}: E \times X \to Y\),
such that ...
\end{definition}

Then, as accustomed with category theory, we need some
universal property. Since \(\mathrm{ev}\) already describes
how to form morphisms \emph{out of} \(E\), our universal
property describes how to create morphisms \emph{into} \(E\):

\begin{definition*}[Continued]
... if there is an object \(S\) with a morphism \(u : S \times X \to Y\),
then there is a unique morphism \(v : S \to E\)
\[\begin{tikzcd}
    S & {S \times X} \\
    \\
    E & {E\times X} && Y
    \arrow["{\mathrm{ev}}", from=3-2, to=3-4]
    \arrow["u", from=1-2, to=3-4]
    \arrow[dashed, from=1-2, to=3-2]
    \arrow["{v}"{description}, dashed, from=1-1, to=3-1]
\end{tikzcd}\]
such that, if the dashed arrow in the triangle is filled
with \(v \times \mathrm{id}\) (which is the Haskell
\texttt{first v = v *** id}), then the diagram commutes.
\end{definition*}

This is basically describing \emph{lambda abstraction}.
Given a function \(u\), we have \(u(s, x) \in Y\), so
we can form the function \(v(s) = \lambda x. u(s, x)\).
\footnote{Note how we use ``pointful'' notation
--- notation involving elements \(x \in X\) etc. ---
to give intuition of ``point-free'' definitions.
In this article it is only a convenient device to describe
\emph{rough feelings} of certain definitions. But in fact,
it can be made rigorous as the \textbf{internal language}
of a topos, where we can freely write expressions like this,
and be confident that they can be traslated back into the
category language.}

The exponential construction creates a functor \((-)^X\). Also,
in Haskell language, the \(\mathsf{fmap}\, f\) instance
of \((-)^X\) is exactly \texttt{(f .)}, the left compositions.

A brilliant insight of exponentials is that they are completely
characterized by \emph{currying}:

\begin{theorem}
There is a natural isomorphism \[\Hom(X \times Y , Z) \cong \Hom(X, Z^Y).\]
In other words, \[(-) \times Y \dashv (-)^{Y}.\]
\end{theorem}

The interested reader shall complete the proof. Next, we continue
on our quest of right inverses. We of course want to
express the \emph{identity morphism} first:

\[\begin{tikzcd}
    {1\times X} \\
    \\
    {X^X \times X} && X
    \arrow["{\mathrm{ev}}"{description}, from=3-1, to=3-3]
    \arrow["{\pi_2}"{description}, from=1-1, to=3-3]
    \arrow[dashed, from=1-1, to=3-1]
\end{tikzcd}\]

Here the dashed line is the unique morphism
\(\mathfrak{id} \times \mathrm{id}\), where
\(\mathfrak{id} : 1 \to X^X\) picks out the identity function
in the object \(X^X\).

Now that we have \(\mathfrak{id}\) as our equipment, consider
this pullback, where \(f : Y \to X\):

\[\begin{tikzcd}
    Z && {Y^X} \\
    \\
    1 && {X^X}
    \arrow["{\mathfrak{id}}"{description}, from=3-1, to=3-3]
    \arrow["{\mathsf{fmap}\,f}"{description}, from=1-3, to=3-3]
    \arrow[from=1-1, to=1-3]
    \arrow[dashed, from=1-1, to=3-1]
    \arrow["\lrcorner"{anchor=center, pos=0.125}, draw=none, from=1-1, to=3-3]
\end{tikzcd}\]

Returning to where we tangented off, the pullback
\(Z\) is, in the category \(\mathsf{Set}\),
the set \(\{g \in Y^X \mid f\circ g = \mathrm{id}\}\).
(Recall that \(g \in Y^X\) means \(g\) is a function \(X \to Y\).)
This captures exactly the right inverses of \(f\).

\subsection*{Fiberwise juggling}

Putting the solution in use, since a fibered space
\(E = \sum_{x:B} p^{-1}(x)\) is defined by a morphism
\(p : E \to B\), we need to find the space of right
inverses of \(p\), which should give the space of dependent products.

\[\begin{tikzcd}
    Z && {E^B} & E \\
    \\
    1 && {B^B} & B
    \arrow["{\mathsf{fmap}\,p}"{description}, from=1-3, to=3-3]
    \arrow["{\mathfrak{id}}"{description}, from=3-1, to=3-3]
    \arrow[dashed, from=1-1, to=3-1]
    \arrow[from=1-1, to=1-3]
    \arrow["p"{description}, from=1-4, to=3-4]
\end{tikzcd}\]

This \(Z\) (considered as a fibered space
\(Z \to 1\)) is then what we sought for.

Now we can generalize from \(\Cat C/1\) to arbitrary
slice categories \(\Cat C/A\).
We are now given a morphism \(f : B \to A\), and we
are supposed to construct a functor \(f_* : \Cat C/B \to \Cat C/A\).
As before, let \(p : E \to B\) be an object of
\(\Cat C/B\). Thinking in \(\mathsf{Set}\)-language,
we should have a ``fiberwise right inverse'' \(p_a^{-1}(x)\),
whose domain is the fiber \(B_a = f^{-1}(a)\) of \(B\) over \(a \in A\).
Its codomain would naturally be \(E_a\), which is a fiber
of \(E\) when considered as a fibered space \((f\circ p) : E \to A\).
Each fiber of the dependent product object \(f_* p\) should look like
\[\prod_{x:B_a} p_a^{-1}(x).\]

The fiberwise right inverse is easy enough to construct (note
that we are still working in \(\mathsf{Set}\)). We just
replace everything in the previous construction.

\[\begin{tikzcd}
    Z_a && {{E_a}^{B_a}} & E_a \\
    \\
    1 && {{B_a}^{B_a}} & B_a
    \arrow["{\mathsf{fmap}\,p_a}"{description}, from=1-3, to=3-3]
    \arrow["{\mathfrak{id}}"{description}, from=3-1, to=3-3]
    \arrow[dashed, from=1-1, to=3-1]
    \arrow[from=1-1, to=1-3]
    \arrow["{p_a}"{description}, from=1-4, to=3-4]
\end{tikzcd}\]

We have the fiberwise constructions ready. How can we ``collect the
fibers'' to create a definition that does not refer to the ``points''
\(a \in A\)? It looks like we are stuck.
Maybe it's time to take a retrospect of what we've achieved.

\section{The True Nature of Slice Categories}

Concepts in category theory are like elephants. You may,
through analogies, theorems, or practical applications,
grasp a feeling of what those concepts are like. But in truth,
these feelings are only describing a part of the elephant.
So let me reveal yet another part, yet another blind man's
description of elephants:

\begin{center}
    \emph{Slice categories descibe local, fiberwise constructs.}
\end{center}

Let's return again to the definition of a fibered product.
\[\begin{tikzcd}
    X \\
    & Z && C \\
    \\
    & B && A
    \arrow["g"{description}, from=2-4, to=4-4]
    \arrow["f"{description}, from=4-2, to=4-4]
    \arrow["h"{description}, from=2-2, to=4-4]
    \arrow[from=2-2, to=4-2]
    \arrow[from=2-2, to=2-4]
    \arrow[dashed, from=1-1, to=2-2]
    \arrow[curve={height=6pt}, from=1-1, to=4-2]
    \arrow[curve={height=-6pt}, from=1-1, to=2-4]
\end{tikzcd}\]
I have added another morphism \(h\), which does not change
the definition since everything commutes in this diagram.
But it brings an interesting change of perspective:
\(h : Z \to A\), regarded as an object in \(\Cat C/A\), is
exactly the usual product of the objects \(f\) and \(g\)!

On second thought this is very natural: Everything in slice
categories needs to respect fibers, i.e. given two fibered spaces
\(B \to A\) and \(C \to A\), any morphisms between them
must map anything in the fiber \(B_a\) over \(a\) to
the fiber \(C_a\). Therefore, the categorical product
of two fibered spaces should also be the fiberwise product.
This immediately generalizes to any construction.

\textsc{Exercise}: Define the notion of fibered coproducts, and
explain why it is the coproduct in the slice category.
Also, explain why the ``fiberwise terminal object'' is exactly
\(\mathrm{id} : A \to A\).

One thing to keep in mind: When we are talking about
the category \(\Cat C/A\), the fibers are considered
to be over \(A\). So when we switch to a different category
\(\Cat C/B\), the spaces are now considered fibered over \(B\).
That's essentially the content of the base change functor:
it changes the base space of the fiber spaces.

Armed with new weapons, we can finally write down the definition
of dependent products:

\[\begin{tikzcd}
    {(Z \to A)} && {(E\stackrel {f\circ p} \to A)^{(B\stackrel f \to A)}} & E \\
    \\
    {(A\stackrel {\mathrm{id}} \to A)} && {(B\stackrel f \to A)^{(B\stackrel f \to A)}} & B
    \arrow["{\mathsf{fmap}\,p}"{description}, from=1-3, to=3-3]
    \arrow["{\mathfrak{id}}"{description}, from=3-1, to=3-3]
    \arrow[dashed, from=1-1, to=3-1]
    \arrow[from=1-1, to=1-3]
    \arrow["p"{description}, from=1-4, to=3-4]
\end{tikzcd}\]

Note that this commutative diagram is entirely
in the slice category \(\Cat C/A\), where each object
are \emph{arrows} in \(\Cat C\). The exponential objects
are also inside the slice category.
This pullback gives a space \(Z \to A\).

According to our guess at the beginning of this section,
we should denote \(Z \to A\) as \(f_* p\). But of course
we need to verify the functorality of this construction.
But it should be clear, since everything used (exponentials,
products and pullbacks) is functorial.

But there is an even more succinct description of
all these: the dependent product functor
is exactly the \textbf{right adjoint} of
the base change functor \(f^*\). The proof is not
hard, although the diagram involved is a bit messy
if you insist on drawing everything in \(\Cat C\)
instead of the slice categories.

\section{Locally Cartesian Closed}

A cartesian closed category is a category where the terminal
object, all binary products and all exponentials exist.
A locally cartesian closed category is a category whose
\emph{slice} categories are all cartesian closed. Let's
unpack the definition and see what this means.

The terminal object in a slice category \(\Cat C/A\) is
exactly \(\mathrm{id} : A \to A\). So it always exists
in slice categories.
A binary product in a slice category, as we have discussed,
is exactly the fibered product, or pullback.
Therefore, a locally cartesian closed category
should have all pullbacks.

What about local exponentials? If there are two
objects \(p : Y \to A\) and \(q : X \to A\), then the local
exponential object \(p^q : E \to A\) should be defined by
the following diagram:

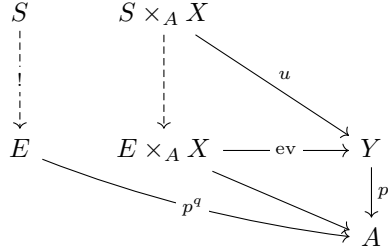
\[\begin{tikzcd}
    S & {S\times_A X} \\
    \\
    E & {E \times_A X} && Y \\
    &&& A
    \arrow["{\mathrm{ev}}"{description}, from=3-2, to=3-4]
    \arrow["p", from=3-4, to=4-4]
    \arrow[from=3-2, to=4-4]
    \arrow["u", from=1-2, to=3-4]
    \arrow["{p^q}"{description}, curve={height=6pt}, from=3-1, to=4-4]
    \arrow[dashed, from=1-2, to=3-2]
    \arrow["{!}"{description}, dashed, from=1-1, to=3-1]
\end{tikzcd}\]

... Well, this looks messy. Let's try the adjoint functor
definition of exponentials: The exponential
functor \((-)^{Y}\) is the right adjoint of the
product functor \((-) \times Y\). So in other words
we should find a right adjoint to the pullback functor
\((-)\times_A Y\). But hey! That looks like the
dependent product functor in the last section.
However, the acute reader may have noticed a discrepancy:
Our dependent product functor is defined as a pullback
of an exponential object. It can't exactly be
the exponential functor, can it?
In fact they have different codomains: Given \(f : C \to A\),
the dependent product functor \(f_* : \Cat C/C \to \Cat C/A\)
is the adjoint of the base change functor
\(f^* : \Cat C/A \to \Cat C/C\). But when we are looking for
the exponential functor, the pullback functor we
want is \((-)\times_A C : \Cat C/A \to \Cat C/A\).
Looking at the diagram for pullbacks we see why:
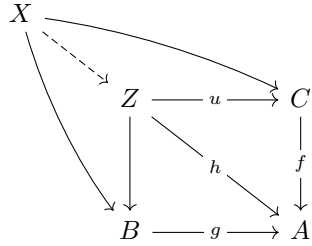
\[\begin{tikzcd}
    X \\
    & Z && C \\
    \\
    & B && A
    \arrow["f"{description}, from=2-4, to=4-4]
    \arrow["g"{description}, from=4-2, to=4-4]
    \arrow["h"{description}, from=2-2, to=4-4]
    \arrow[from=2-2, to=4-2]
    \arrow["u"{description}, from=2-2, to=2-4]
    \arrow[dashed, from=1-1, to=2-2]
    \arrow[curve={height=6pt}, from=1-1, to=4-2]
    \arrow[curve={height=-6pt}, from=1-1, to=2-4]
\end{tikzcd}\]
The functor \((-)\times_A C : \Cat C/A \to \Cat C/A\)
sends \(g\) to \(h\), while the functor \(f^*\) sends
\(g\) to \(u\). Since \(h = f\circ u\), you can see
that the functor \((-)\times_A C\) is the composition
of two functors \(f_! f^*\).

Now we can save a tremendous amount of work
with this theorem:

\begin{theorem}
Given two adjoint pairs:
\[\begin{tikzcd}
    {\Cat C} && {\Cat D} && {\Cat E}
    \arrow[""{name=0, anchor=center, inner sep=0}, "{F_1}"{description}, curve={height=-12pt}, rightarrow, from=1-1, to=1-3]
    \arrow[""{name=1, anchor=center, inner sep=0}, "{F_2}"{description}, curve={height=-12pt}, rightarrow, from=1-3, to=1-5]
    \arrow[""{name=2, anchor=center, inner sep=0}, "{G_1}"{description}, curve={height=-12pt}, rightarrow, from=1-3, to=1-1]
    \arrow[""{name=3, anchor=center, inner sep=0}, "{G_2}"{description}, curve={height=-12pt}, rightarrow, from=1-5, to=1-3]
    \arrow["\dashv"{anchor=center, rotate=-90}, draw=none, from=0, to=2]
    \arrow["\dashv"{anchor=center, rotate=-90}, draw=none, from=1, to=3]
\end{tikzcd}\]
The composition also forms an adjunction
\[F_2F_1 \dashv G_1G_2.\]
\end{theorem}
\begin{proof}%
\[\Hom_{\Cat E}(F_2F_1X, Y) \cong
\Hom_{\Cat D}(F_1X, G_2Y) \cong
\Hom_{\Cat C}(X, G_1G_2 Y).\qedhere\]
\end{proof}

\[\begin{tikzcd}
    {\Cat C / A} &&& {\Cat C/C} &&& {\Cat C/A}
    \arrow["{(-)\times_A Y}"{description}, curve={height=30pt}, from=1-1, to=1-7]
    \arrow[""{name=0, anchor=center, inner sep=0}, "{f_!}"{description}, curve={height=12pt}, from=1-4, to=1-7]
    \arrow["{?}"{description}, curve={height=30pt}, from=1-7, to=1-1]
    \arrow[""{name=1, anchor=center, inner sep=0}, "{f^*}"{description}, curve={height=12pt}, from=1-7, to=1-4]
    \arrow[""{name=2, anchor=center, inner sep=0}, "{f^*}"{description}, curve={height=12pt}, from=1-1, to=1-4]
    \arrow[""{name=3, anchor=center, inner sep=0}, "{f_*}"{description}, curve={height=12pt}, from=1-4, to=1-1]
    \arrow["\dashv"{anchor=center, rotate=90}, draw=none, from=0, to=1]
    \arrow["\dashv"{anchor=center, rotate=90}, draw=none, from=2, to=3]
\end{tikzcd}\]

With this diagram it is crystal clear that the
fibered exponential fits exactly in the position of the
question mark.

In fact, the condition that the
dependent sum functor \(f_!\) (which exists in
every category) has a chain of three adjoints
\[f_! \dashv f^* \dashv f_*, \]
is equivalent to
the condition that the category is locally cartesian closed.
The backward implication is precisely what we proved
in the last section. As for the forward implication,
it is proved by our discussion in the previous few paragraphs.

\section{Prospects}

This introduction has gotten way too lengthy. But I shall
point out several direction to proceed before I end.

Cartesian closed category, as can be seen in the definition,
serves as the semantics of simply typed lambda calculus.
You might not be able to figure out the details at once,
but you should see that there is a probable connection here.
On the other hand, \emph{locally} cartesian closed categories
are central to the semantic interpretation of \emph{dependent
types}. Type dependency is, fundamentally, expressing
fiber spaces; working with dependent types amounts to
making fiberwise comstructions. The classical reference for this
is \cite{seely}.

Although I did not mention any topology in the text, fiber
spaces ultimately came from topology. And it is the fact that
there is a notion of ``neighbourhoodness'' between fibers
that makes them important --- otherwise they are just random
sets.

Going further in this direction,
the adjunction \(f^* \dashv f_*\) is called a \textbf{geometric
morphism} in the language of topos. If it has further adjoints,
it becomes ``smoother'' in the geometric sense. This plays
the central role in topos theory. More can be read at \cite{sketches}.

\end{document}